\documentclass[3p, sort&compress]{elsarticle}

\makeatletter
\def\ps@pprintTitle{%
 \let\@oddhead\@empty
 \let\@evenhead\@empty
 \def\@oddfoot{\centerline{\thepage}}%
 \let\@evenfoot\@oddfoot}
\makeatother

\journal{Journal of Computational Physics}

\usepackage{amsmath}
\usepackage{amssymb}
\usepackage{amsthm}
\usepackage{mathrsfs}
\usepackage{subcaption}
\usepackage{siunitx}
\usepackage{color}
\usepackage[linesnumbered,ruled]{algorithm2e}

\usepackage[mathlines]{lineno}
\usepackage{hyperref}
\modulolinenumbers[5]

\newcommand{\R}{\mathbb{R}}

\newcommand{\E}{\mathcal{E}}
\DeclareMathAlphabet\mathbfcal{OMS}{cmsy}{b}{n}

\newcommand{\HO}[1]{H^{#1}(\Omega)}
\newcommand{\HpO}[1]{H^{#1}(\partial\Omega)}
\newcommand{\LO}[1]{L^{#1}(\Omega)}
\newcommand{\Div}{\nabla\cdot}
\newcommand{\norm}[1]{\left\lVert#1\right\rVert}

\DeclareMathOperator{\ran}{Ran}
\DeclareMathOperator*{\argmin}{arg\,min}

\biboptions{}
\begin{document}

\begin{frontmatter}

\title{Levenberg-Marquardt algorithm for acousto-electric tomography \\based on the complete electrode model}

\author[Changyou]{Changyou Li}
\author[DTU]{Mirza Karamehmedovi\'c}
\author[Katja]{Ekaterina Sherina}
\author[DTU]{Kim Knudsen}
\address[Changyou]{School of Information and Electronics, Northwestern Polytechnical University}
\address[DTU]{Department of Applied Mathematics and Computer Science, Technical University of Denmark}
\address[Katja]{Department of Mathematics, University of Vienna}

\begin{abstract}
  The inverse problem in Acousto-Electric tomography concerns the reconstruction of the electric conductivity in a body from knowledge
of the power density function in the interior of the body.  This
interior power density results from currents prescribed at boundary
electrodes, and it can be obtained through electro-static boundary
measurements together with auxiliary acoustic probing. Previous works on Acousto-Electric tomography used the continuum model for the electrostatic boundary conditions; however, from  
Electrical Impedance Tomography is it known that the complete electrode model is much more realistic and accurate. 

  In this paper the inverse problem of Acousto-Electric
tomography is posed using the (smoothened) complete electrode model, and a reconstruction method based on the
Levenberg-Marquardt iteration is formulated in appropriate function
spaces. This results in a system of partial differential equations to
be solved in each step. To increase the computational efficiency
and stability, a strategy based on both the complete electrode model
and the continuum model is proposed.

The method is implemented numerically for a two
dimensional scenario, and the algorithm is tested on two different
numerical phantoms, a heart and lung model and a human brain
model. Several numerical experiments are carried out confirming the
feasibility, accuracy and stability of the developed method.

\end{abstract}

\date{}


\end{frontmatter}


\section{Introduction}
Electrical Impedance Tomography (EIT) is an emerging technology for
obtaining the internal conductivity of a physical body from
boundary measurements of currents or voltages on the surface of the
body \cite{Holder2005,Yorkey1990,Calderon1980}. The reconstruction problem in EIT is an
ill-posed problem due to the fact that boundary measurements show little sensitivity to (even large) changes of interior conductivity distribution
\cite{Ammari2008}. Intensive research exists on this topic
\cite{Uhlmann2009,Brown2009}; many regularization methods have
been proposed to overcome the ill-posedness and to improve the imaging
quality \cite{Han1999,Hsiao2001,Chung2005,Knudsen2009}. 

More recently it has been suggested to augment the measurement setup
in EIT with an ultra-sonic device thus yielding the hybrid imaging
method known as Acousto-Electric tomography (AET)
\cite{Ammaribook2008}.  The resulting modality has been investigated
theoretically and numerically, and AET seems to have the potential to
dramatically increase the contrast, resolution, and stability of the
conductivity reconstruction \cite{Bal2012,Kuchment2012}.

 The idea of AET it is to conduct a usual EIT experiment while a known focused ultrasonic wave propagates through the object. The high intensity of the acoustic pressure  will create a small local deformation in the physical body and thus of the electrical conductivity due to the acousto-electric effect \cite{Geselovitz1971,Jossinet1999,Jossinet2005}. 

 The physical body to be imaged is modeled as an open, bounded and
 smooth domain $\Omega\subset\R^n$ with $n\ge 2,$ and the electric
 conductivity in $\Omega$ is described by the smooth function
 $\sigma\geq c > 0.$ When an electric field is applied to the boundary
 $\partial \Omega$ of the body, a voltage potential $u(\sigma)$ is
 generated inside $\Omega.$ With the assumption that measurements are
 carried out in low temporal frequency and that $\Omega$ contains no
 interior sources or sinks of charges the governing equation is the the
 generalized Laplace equation
 \begin{equation}
   \nabla\cdot\sigma\nabla u = 0 \quad \text{in }\Omega\label{eq:epotential}
 \end{equation}
 subject to suitable boundary conditions.

 The acoustic wave perturbs the electric conductivity and consequently  the interior voltage potential. From standard EIT measurements recorded while the wave propagates through the body,  the interior power density 
\begin{equation}
  \label{eq:pdensity}
  \E(\sigma) = \sigma|\nabla u(\sigma)|^2
\end{equation}
can be found  \cite{Ammari2008,Kuchment2012,Jensen2019}. The inverse problem in AET is thus to find $\sigma$ from knowledge of one or several $\E(\sigma)$ corresponding to different boundary conditions.

Several methods have been developed in the literature for reconstructing $\sigma$ from $\E(\sigma)$ \cite{Ammari2008,Ammaribook2008,Kuchment2012,Bal2013,Hoffmann2014,Song2016,Song2017}. All these methods are based on the so-called continuum model, i.e.\ the boundary conditions for \eqref{eq:epotential} are given by continuous fileds as either Dirichlet or Neumann conditions without explicit electrode modelling. Given  measurements $\E^\delta$ of a true power densities $\E$ perturbed by noise of magnitude $\delta,$ a reasonable approach to the reconstruction problem is through optimization
\begin{equation}
  \label{eq:mini}
  \min_\sigma\|\E^\delta - \E(\sigma)\|_{L^2(\Omega)}^2.
\end{equation}
Since the mapping from $\sigma$ to $\E(\sigma)$ is non-linear,
\eqref{eq:mini} is a non-linear least-squares problem in $L^2(\Omega)$
and for that reason and iterative approach such as the
Levenberg-Marquardt method is suitable \cite{Bal2012}.  See
\cite{Adesokan2018,Jensen2018,roy2018a,gupta2019a} for alternative
optimization approaches, and \cite{Hubmer2018} for a related analysis of the
limited boundary data problem.

The Complete Electrode Model (CEM) is a practical model for EIT that explicitly models electrodes; CEM can simulate EIT measurements with much greater accuracy than the continuum models \cite{Somersalo1992}. In the model, $L$ electrodes are attached on boundary $\partial\Omega$. A known total current $I_l$ injected through the $l$-th electrode $e_l$ is given as
\begin{equation}
  \label{eq:bc1}
  \int_{e_l}\sigma\frac{\partial u}{\partial \nu} = I_l, \quad l=1,\ldots,L.
\end{equation}
Here, $\nu$ is the outward unit normal vector to the boundary $\partial\Omega,$ $e_l$ is the  $l$-th electrode, and ${\partial}/{\partial \nu}$ indicates the derivative of $u$ in the direction of the  outward unit normal vector $\nu$. Since there is no current flowing out through boundary regions without electrodes, one has
\begin{equation}
  \label{eq:bc2}
  \sigma\frac{\partial u}{\partial \nu} =0 \quad \text{on } \partial\Omega\backslash\bigcup^{L}_{l=1}e_l.
\end{equation}
On the electrode $e_l$ the electric potential $U_l$ is assumed to be constant (but unknown). This boundary  potential consist of a part due to the interior potential $u$ and a part due to the electrode contact, and it is comprised in the model 
\begin{equation}
  \label{eq:bc3}
  u + z_l\sigma\frac{\partial u}{\partial \nu} =  U_l \quad\text{on } e_l,
\end{equation}
where $z_l$ denotes the so-called contact impedance assumed to be constant on the $l$-th electrode. The partial differential equation \eqref{eq:epotential} with boundary conditions \eqref{eq:bc1}-\eqref{eq:bc3} gives the CEM. To ensure existence and uniqueness of the solution, this model also needs to include the law of charge conservation
\begin{equation}
  \label{eq:bc4}
  \sum^{L}_{l=1}I_l = 0
\end{equation}
and to determine the potential's grounding  by
\begin{equation}
  \label{eq:bc5}
  \sum^{L}_{l=1}U_l = 0.
\end{equation}
The CEM problem \eqref{eq:epotential} with \eqref{eq:bc1}-\eqref{eq:bc5} has a unique solution $(u,U),$ $U = (U_1,U_2,\ldots,U_L),$ with  $u \in \HO{2-\epsilon}$ for any $\epsilon > 0$  \cite{Hyvonen2017}, however the mixed  boundary conditions \eqref{eq:bc2}-\eqref{eq:bc3} allow singularities near the edges of the electrodes. Theoretically and computationally such singularties are challenging. In order to increase the regularity of $u,$ the electrode conductance $\zeta_l = 1/z_l$ is introduced in \cite{Hyvonen2017} as a smoooth function, thus \eqref{eq:bc3} is replaced by
   \begin{align}
       \sigma\frac{\partial u}{\partial \nu} &= \zeta(U_l-u) \quad \text{on } e_l. \label{eq:scem_model}
  \end{align}
The PDE problem  \eqref{eq:epotential} with \eqref{eq:bc1}-\eqref{eq:bc2} and \eqref{eq:bc4}-\eqref{eq:scem_model} is called the smoothened CEM (SCEM).  When  $\zeta_l \in \HpO{t}$ for some $t > (n-1)/2,$ the potential $u(\sigma) \in \HO{t+\frac{3}{2}}$ \cite{Hyvonen2017}. In particular for $s=t+1/2>{n}/{2}$,  $u(\sigma)\in \HO{s+1}$ and $\E(\sigma) \in \HO{s},$ since the Sobolev space $H^s(\Omega)$ is a Banach algebra when $s > n/2.$

The aim of this paper is to develop an iterative method for reconstructing $\sigma$ from $\E^\delta.$ Inspired by \cite{Bal2012} the method will be based on the Levenberg-Marquardt method for solving the least squares problem \eqref{eq:mini}. In contrast to \cite{Bal2012}, who relied on a continuum model, the main novelty here is the use the SCEM to accurately and stably model electrodes and the electric current in the forward problem.
In addition, a reconstruction strategy based on a combined use of both CEM and continuum model with Dirichlet boundary condition (DCM) is  proposed to increase the computational efficiency.

The outline of the paper is as follows. In Section \ref{sec:2}, the Levenberg-Marquardt method is briefly introduced. The non-linear problem is linearized, and the adjoint problem is formulated. In Section \ref{sec:numeri_imple}, the iterative reconstruction method is developed based on CEM and Levenberg-Marquardt iteration. A linear system is build that calculates the updating step for each iteration. The algorithm for increasing efficiency by exploiting DCM is also introduced in this section. These algorithms are implemented and applied to reconstruct the conductivity distribution of several phantoms in Section \ref{sec:3}, numerical performances of different algorithms are discussed in detail. The conclusion of the presented work is given in Section \ref{sec:4}.

\section{Reconstruction algorithm}\label{sec:2}
In this section we will first recap the Levenberg-Marquardt Algorithm (LMA) for solving (non-linear) optimization problems. Then we will for the particular problem in AET with CEM derive the necessary ingredients, that is the Fr\'echet derivative of $\E$ and it's adjoint.
\subsection{The Levenberg-Marquardt Algorithm}
Let $F: X\rightarrow Y$ be a (possibly non-linear) operator between Hilbert spaces $X$ and $Y.$ For some $y^\delta=y + \delta \in Y$ with $y \in \ran(F)$ the problem is to solve at least approximately the equation $F(\sigma) \approx y^\delta,$ and often the minimization problem
\begin{align}\label{eq:mini_NL}
  \argmin_{\sigma\in X}\|F(\sigma) - y^\delta\|^2_Y
\end{align}
is considered. If $F$ is (Fr\'echet) differentiable, the linear approximation 
\begin{equation*}
   F(\sigma) - F(\sigma_{k}) \approx F'(\sigma_k)(\sigma-\sigma_k) 
\end{equation*}
yields the iterative scheme
\begin{equation*}
  \sigma_{k+1} = \argmin_{\sigma} \norm{y^\delta - F(\sigma_k) - F'(\sigma_k)(\sigma-\sigma_k)}^2_Y
\end{equation*}
solved by
\begin{align*}
 \sigma_{k+1} = \sigma_k + F'(\sigma_k)^{-1}(y^\delta - F(\sigma_k)).
\end{align*}
This is applicable only when $F'(\sigma_k)$ is left-invertible. In general one can instead minimize the following Tikhonov functional
\begin{equation}
  \label{eq:minipro}
  \sigma_{k+1} = \argmin_{\sigma}{\norm{y^\delta - F(\sigma_k) - F'(\sigma_k)(\sigma-\sigma_k)}^2_Y + \alpha_k\norm{\sigma - \sigma_k}^2_X},
\end{equation}
where $\alpha_k>0$ is the regularization parameter; this minimization problem is solved by  
\begin{equation}
  \label{eq:lm_iter}
  \sigma_{k+1}=\sigma_k+(F'(\sigma_k)^*F'(\sigma_k) + \alpha_kId)^{-1}F'(\sigma_k)^*(y^{\delta} - F(\sigma_k))
\end{equation}
known as the Levenberg-Marquardt Algorithm (LMA).
Here, $F'(\sigma_k)^*$ is the adjoint of $F'(\sigma_k)$ and $Id$ is identity operator. When the operator $F$ satisfies a certain nonlinarity condition, a proper choice of the parameter $\alpha_k$ and the initial guess $\sigma_0$ sufficiently close to the desired solution, the LMA converges to a solution $\sigma^\delta$ of \eqref{eq:mini_NL} \cite{Hanke1997,Kaltenbacher2008}.

 The LMA  can be thought of as a combination of steepest descent and Gauss-Newton method. When the current solution is far from the correct one, a large value is assigned to $\alpha_k>0$, and LMA behaves like a steepest-descent method which converges slowly. When the current solution is close to the correct solution, a small $\alpha_k>0$ is used, and LMA behaves like a Gauss-Newton method, which has faster convergence.


\subsection{The Fr\'echet derivative $\E'$ and its adjoint}
In the present work we consider the operator  $F(\sigma) = \E(\sigma)$ defined in \eqref{eq:pdensity}. In order to apply \eqref{eq:lm_iter} we therefore need to calculate the Fr\'echet derivative $\E'(\sigma)$  of the power density operator at $\sigma$ and its adjoint $\E'(\sigma)^\ast.$ 

We start with the Fr\'echet derivative  of $u(\sigma),$ i.e.\ the operator $u'(\sigma): H^s(\Omega) \rightarrow H^{s+1}(\Omega).$ The Fr\'echet derivative can be calculated (see the general approach in \cite{Lechleiter2008a}) in the following way: 
For a given $\tau \in H^s(\Omega)$ with compact support inside $\Omega,$ the difference $u(\sigma +\tau) - u(\sigma)$ is approximated by a function $\xi = u'(\sigma)\tau$ that is linear in $\tau.$ Indeed, define $(\xi,\Xi) \in H^{s+1}(\Omega) \times \R^n$  as the solution to the modified CEM problem
\begin{subequations}
  \label{eq:final_frechet}
  \begin{align}
  \label{eq:spoisson_xi}
  \nabla\cdot \sigma\nabla \xi + \nabla\cdot(\tau\nabla u(\sigma)) &=0 \quad\text{in }\Omega,\\
  \label{eq:sbc1_xi}
  \sigma\frac{\partial \xi}{\partial \nu} &= \zeta(\Xi_l-\xi) \quad\text{on } e_l,\\
  \label{eq:sbc2_xi}
    \int_{e_l}\sigma\frac{\partial \xi}{\partial \nu} &= 0 \quad\text{on } e_l,\\
    \sigma\frac{\partial \xi}{\partial \nu} &= 0 \quad\text{off } e_l.
  \end{align}
\end{subequations}
With the grounding  $\sum \Xi_l = 0,$ \eqref{eq:final_frechet} has a unique weak solution $\xi \in \HO{s+1},\; \Xi\in \R^n,$ which obviously is linear with respect to $\tau$. Moreover, 
\begin{align*}
  \|u(\sigma +\tau) - u(\sigma) - \xi \|_{H^{s+1}(\Omega)} \leq C\|\tau\|_{H^s(\Omega)}
\end{align*}
showing that $\xi = u'(\sigma)\tau$ is indeed the Fr\'echet derivative of $u$ at $\sigma.$
The Fr\'echet derivative $\E'(\sigma)\tau$ of $\E(\sigma)$ in the direction $\tau$ is obtained as in \cite{Bal2013} now given by  
\begin{equation}
  \label{eq_eptau}
    \mathcal{E}'(\sigma)\tau = \tau|\nabla u(\sigma)|^2 +2\sigma\nabla u(\sigma)\cdot \nabla \xi.
\end{equation}


We now compute the adjoint $(\E'(\sigma))^\ast$ first as an operator in $L^2(\Omega):$ Consider for some $z \in L^2(\Omega)$
\begin{align}\label{eq:adj1}
  \langle z, \E'(\sigma)\tau \rangle_{L^2(\Omega)} = \int z\left( \tau |\nabla u(\sigma)|^2 + 2\sigma \nabla u(\sigma) \cdot \nabla \xi\right) dx.
\end{align}
We focus on the second term in the integral (the first term is self-adjoint). Introduce (based on experience) the auxilary pair $v(z) \in H^1(\Omega)$ and $V \in \R^n$ defined by the weak PDE form (for all $w,W$)
\begin{align}\label{eq:vV}
    \int_\Omega \sigma \nabla v(z) \cdot \nabla w \; dx + \sum_l \int_{\partial \Omega} \zeta (V_l(z) - v(z))(W_l-w)\; dS = \int_\Omega 2 \sigma z  \nabla u(\sigma) \cdot \nabla w \; dx.
\end{align}
The strong form reads
\begin{align*}  
    \nabla \cdot \sigma \nabla v -2 \nabla \cdot \sigma z  \nabla u(\sigma) &= 0, \\
\sigma \frac{ \partial v }{\partial \nu} &= \zeta (V_l-v) \; \text{on} \; e_l,\\ 
\int_{e_l}\sigma \frac{\partial v}{\partial \nu} dS &= 0,\\ 
\sigma \frac{\partial v}{\partial \nu} &= 0 \; \text{off} \; e_l.
\end{align*}

Inserting the pair $(w,W) = (\xi,\Xi)$ in \eqref{eq:vV} we calculate the latter term in the right hand side of \eqref{eq:adj1} 
\begin{align*}
  \int 2z \sigma \nabla u(\sigma) \cdot \nabla \xi dx &=   \int_\Omega \sigma \nabla v(z) \cdot \xi dx + \sum_l \int_{\partial \Omega} \zeta (V_l(z) - v(z))(\Xi_l - \xi)\; dS\\
&=  - \int_\Omega \tau \nabla u(\sigma)  \cdot \nabla v(z)  \; dx,
\end{align*}
where the last equality follows from the weak form of \eqref{eq:final_frechet}. Thus we find
\begin{align*}
\langle (\E'(\sigma))^* z,\tau\rangle_{L^2(\Omega)} & =  \langle z, \E'(\sigma)\tau \rangle_{L^2(\Omega)}\\
&=  \langle  |\nabla u(\sigma)|^2z - \nabla u(\sigma)  \cdot \nabla v(z), \tau \rangle_{L^2(\Omega)},
\end{align*}
that is
\begin{align}\label{eq:tehz}
  (\E'(\sigma))^* z = |\nabla u(\sigma)|^2z - \nabla u(\sigma)  \cdot \nabla v(z).
\end{align}

To get to the adjoint in a higher order Sobolev spaces $H^s(\Omega)$ for the chosen $s>n/2$ (e.g.\  $s = 2$ in 2D) we lift the operator to higher order spaces, i.e.\ we solve for $x \in H^s(\Omega)$ the equation
\begin{align*}
  \langle x,\tau \rangle_{H^s(\Omega)} =  \langle (\E'(\sigma))^* z,\tau\rangle_{L^2(\Omega)}.
\end{align*}
Using the embedding operator $B: \HO{s}\rightarrow \LO{2}$ we can write  $x = B^*\E'(\sigma)^* \in \HO{s},$ with $B^*$ denoting the Banach space adjoint of $B$ and $(\E'(\sigma))^*z$ from \eqref{eq:tehz}. This is a fourth order PDE problem when $s=2.$

\section{Iterative reconstruction algorithm based on LMA}
\label{sec:numeri_imple}
According to the Levenberg-Marquardt iteration given in \eqref{eq:lm_iter}, the formulation for calculating the $k$-th updating step $\tau_k$ for the presented problem is explicitly given as
\begin{equation}
  \label{eq:tauk}
  (\E'(\sigma_k)^*\E'(\sigma_k) + \alpha_kId)\tau_k = \E'(\sigma_k)^*(\E^{\delta} - \E(\sigma_k)),
\end{equation}
where $\E'(\sigma_k)^*\E'(\sigma_k)\tau = B^*M\tau$ and
\begin{equation*}
  M\tau \\= |\nabla u|^2(\tau|\nabla u|^2 +2\sigma\nabla u\cdot \nabla u'(\sigma)\tau)+ 2\nabla u \cdot \nabla v (\tau|\nabla u|^2) + 4\nabla u \cdot \nabla v(\sigma\nabla u\cdot \nabla u'(\sigma)\tau).
\end{equation*}
Here, $M\tau$ is easily obtained with \eqref{eq_eptau} and \eqref{eq:tehz}. After computing $\tau_k$ from \eqref{eq:tauk}, the conductivity map $\sigma_k$ obtained from the $k$-th iteration is updated by $\sigma_{k+1} = \sigma_k+\tau_k$ for a new iteration. 
All PDEs are coupled and collected into the PDE system
\begin{subequations}
  \label{eq:final_lma_final}
\begin{align}
  \label{eq:final}
  \phi + \alpha_k\tau_k &=y, \qquad\text{in } \Omega,\allowdisplaybreaks\\
  \Delta \phi - \chi &= 0,\qquad\text{on}\ \partial\Omega, \allowdisplaybreaks\\
  \beta^2\Delta\chi + \phi - \gamma &= 0, \qquad\text{on}\ \partial\Omega,\allowdisplaybreaks\\
  \partial \phi /\partial \nu &= 0, \qquad\text{on}\ \partial\Omega\allowdisplaybreaks\\
  \partial \chi /\partial \nu &= 0, \qquad\text{on}\ \partial\Omega\allowdisplaybreaks\\  
  \gamma - |\nabla u|^2 \tau|\nabla u|^2 +2\sigma\nabla u\nabla \xi+ 2\nabla u \nabla \zeta + 4\nabla u \cdot \nabla \kappa &=0\allowdisplaybreaks, \qquad\text{in } \Omega,\\
  \label{eq:final_eq1}
  \nabla\cdot\sigma\nabla \xi + \nabla\cdot(\tau\nabla u(\sigma))&=0, \qquad\text{in } \Omega, \allowdisplaybreaks\\
  \label{eq:final_bc1}
  \sigma\frac{\partial \xi}{\partial \nu} &= \zeta(\Xi-\xi), \qquad\text{on } e_l,\allowdisplaybreaks\\
  \label{eq:final_bc2}
  \int_{e_l}\sigma\frac{\partial \xi}{\partial \nu} &= 0, \qquad\text{ on } e_l,\allowdisplaybreaks\\
  \nabla\cdot\sigma\nabla\rho+\Div(\tau|\nabla u|^2\sigma\nabla u(\sigma))&=0, \ \qquad\text{in } \Omega, \allowdisplaybreaks\\
  \sigma\frac{\partial \rho}{\partial \nu} &= \zeta(\varrho_l-\rho) \qquad\text{on } e_l, \allowdisplaybreaks\\
  \int_{e_l}\sigma\frac{\partial \rho}{\partial \nu} &= 0 \qquad\text{ on } e_l,  \allowdisplaybreaks\\
  \sigma\frac{\partial \rho}{\partial \nu} &= 0 \qquad\text{ off } e_l,\allowdisplaybreaks\\
  \nabla\cdot\sigma\nabla\kappa + \Div((\sigma\nabla u\nabla \xi)\cdot\sigma\nabla u(\sigma))&=0, \qquad\text{in } \Omega, \allowdisplaybreaks\\
  \label{eq:final_bck}
  \sigma\frac{\partial \kappa}{\partial \nu} &= \zeta(\varkappa_l-\kappa) \quad\text{on } e_l, \allowdisplaybreaks\\
  \int_{e_l}\sigma\frac{\partial \kappa}{\partial \nu} &= 0, \qquad\text{on } e_l,  \allowdisplaybreaks\\
  \sigma\frac{\partial \kappa}{\partial \nu} &= 0 \ \qquad\text{off } e_l
\end{align}
\end{subequations}
%
with $y=\E'(\sigma_k)^*(\E^{\delta} - \E(\sigma_k))$, $\phi = B^*M\tau$, $\gamma = M\tau$, $\xi = u'(\sigma)\tau$, $\rho = v(\tau|\nabla u|^2)$, and $\kappa = v(\sigma\nabla u\cdot\nabla u'(\sigma)\tau)$. With the number of measurements $M>1$, the system is formulated with $y=\mathbfcal{E}(\sigma_k)^*(\mathbfcal{E}^{\delta} - \mathbfcal{E}(\sigma_k))$ and $$\gamma = \sum^M_{m=1}\left[|\nabla u_m|^2\tau|\nabla u_m|^2 +2\sigma\nabla u_m\cdot \nabla \xi_m+ 2\nabla u_m \cdot \nabla \zeta_m + 4\nabla u_m \cdot \nabla \kappa_m\right].$$ Since equations \eqref{eq:final_eq1}-\eqref{eq:final_bck} need to be solved for each measurement, one additional measurement will need 3 additional partial differential equations. 

It may be possible to simplify the above PDE system if the interior potential near the measurement boundary and the measured potential both turn out to converge significantly faster (with the number of iterations of the solution procedure) than the conductivity estimate $\sigma$. In this case, we would expect $\xi\approx0$ and $\Xi\approx0$ to hold early in the iteration. From~\eqref{eq:final_bc1}, we would then get that the change in the current at any $e_l$ resulting from a change in $\sigma$ is approximately zero after only a few iterations, and conditions~\eqref{eq:final_bc1} and~\eqref{eq:final_bc2} might justifiably be substituted with the much simpler Dirichlet boundary condition $\xi=0$ on $\partial \Omega.$

The iterative reconstruction method based on the above system is here named as LM-SCEM which is demonstrated in Algorithm \ref{algo_lm_cem} for a single measurement. The measured data of {$\mathcal{E}^\delta$} is simulated with SCEM in this paper. The relative error is given by $\eta =\|\sigma_{\rm t} - \sigma_{\rm r}\|_{L^2(\Omega)}/\|\sigma_{\rm t}\|_{L^2(\Omega)},$ where $\sigma_{\rm t}$ and $\sigma_{\rm r}$ denotes the true and reconstructed conductivities. The parameter $\alpha_k$ should theoretically be updated according to the value of $\tau_k$. If $\sigma_k + \tau_k$ leads to a reduction of the relative error in $\sigma_k$, $\alpha_k$ is decreased and $\tau_k$ is accepted. Otherwise, $\tau_k$ is discarded and $\alpha_k$ is increased. Since $\eta$ can not be determined in practice, a relatively large value is asisgned to $\alpha_0$ in this paper, and $\alpha_k$ is slowly decreased to ensure convergence. The iteration is stopped when the $L^2$-norm of $\tau_k$ is smaller than a given value or the maximum number of iteration is achieved. 

\begin{algorithm}[H]
  \SetAlgoLined
  \KwData{The measured power density $\mathcal{E}^\delta$ and an initial guess $\sigma_0$}
  \KwResult{The reconstructed conductivity map $\sigma_r$ with a relative error $\eta$ }

  $N$: the maximum number of iterations\;
  $\sigma_k\leftarrow\sigma_0$\;
  $\alpha_k \leftarrow \alpha_0$\;
  $num \leftarrow 1$\;
  $norm \leftarrow 1$\;
  \While{$norm > \delta$ and $num < N$}{
    Update $u_k$ from $\sigma_k$ with CEM\;
    Compute $\mathcal{E}(\sigma_k)$ from $u_k$\;
    Compute $y = \E'(\sigma_k)^*(\E^{\delta} - \E(\sigma_k))$ 
    Compute $\tau_k$ with the linear system defined by \eqref{eq:final_lma_final}\;
    $\sigma_k\leftarrow\sigma_k+\tau_k$\;
    $norm \leftarrow \|\tau_k\|_{L^2}$\;
    $num \leftarrow num+1$\;
    Update $\alpha_k$\;
  }
  
  \caption{The LM-SCEM algorithm for reconstructing the conductivity map from single measurement of power density.}
  \label{algo_lm_cem}
\end{algorithm}

\subsection{The LM-DCM method and the mixed reconstruction algorithm}
If DCM is considered instead of SCEM, the system for reconstructing $\sigma$ can be built by a similar calculation \cite{Bal2013}. The resulted system for computing $\tau_k$ remains the same as the one for LM-SCEM except that the boundary conditions \eqref{eq:final_bc1} and \eqref{eq:final_bc2} need to be replaced with $\xi = 0$ on $\partial \Omega$. The computation of $\tau_k$ in LM-SCEM is actually more expensive because of the additional unknowns in $\Xi$. To obtain a good reconstruction accuracy of LM-SCEM, multiple measurements are usually considered, but the computational efficiency will decrease with the increasing number of measurements. Our investigation on the convergence of boundary potential and the conductivity shows that the boundary potential converges much faster, examples are given in numerical experiments. This is mainly because EIT measurement is not very sensitive to the internal change of the conductivity distribution \cite{Ammari2008}. This property renders EIT an ill-posed problem, but it will be taken as the foundation here to build a faster reconstruction approach by mixing LM-SCEM and LM-DCM, which is abbreviated by LM-SCEM-DCM.

This mixed reconstruction method is illustrated in Algorithm \ref{algo_lm_cem_dcm} for a single measurement. Here, $\mathbfcal{U}=[U_1, U_2,\ldots, U_L]$, which is a vector composed of the voltages on the electrodes. The true values $\mathbfcal{U}^t$ can be measured, which is produced when simulating the power density with CEM, therefore no additional computation is required. The information of  $\mathbfcal{U}^t$ is here exploited to define a stopping criteria for reconstructing the boundary potential $u_b$ with LM-SCEM. With a relative error given by $\eta^b=\frac{\|\mathbfcal{U}^t-\mathbfcal{U}^k\|_{l^2}}{\|\mathbfcal{U}^t\|_{l^2}}$, the value of $\eta^b$ is checked in each iteration, and LM-SCEM is terminated when an expected relative error $\eta^b_0$ is achieved. Since the regularity of potentials computed from CEM is not so good on $\partial \Omega$, a smaller region $\Omega'=\{x\ |\ x\in\Omega,\ dist(x, \partial\Omega)>d\}$ is defined with a small $d$ which can ``smooth'' out the possible irregularity of $u$ close to $\partial\Omega$. The boundary potential $u_b$ is here defined on $\partial\Omega'$ and used for the reconstruction with LM-DCM. The power density in $\Omega'$ is firstly reconstructed from $u_b$. This can be achieved with the method introduced in \cite{Ammari2008} but will be simulated with DCM here. Because there are no additional unknowns $\Xi$ in LM-DCM, the computation will be more efficient, especially for the computation with multiple measurements. Meanwhile, the conductivity map produced by LM-SCEM is used as the initial guess for LM-DCM. This good initial guess will also help LM-DCM converge faster. Therefore, the mixed reconstruction method can provide a practical and efficient computational model for AET.
\begin{algorithm}[t]
  \SetAlgoLined
  \KwData{The measured power density $\mathcal{E}^\delta$ and the voltage vector $\mathbfcal{U}^t$ on electrodes. An initial guess $\sigma_0$ and an expected relative error $\eta^b_0$ for $\mathbfcal{U}^k$ }
  \KwResult{The reconstructed conductivity map $\sigma_r$ with a relative error $\eta$ }

  $N_{c}$: the maximum number of iteration for LM-SCEM\;
  $N_{d}$: the maximum number of iteration for LM-DCM\;
  $\sigma_k\leftarrow\sigma_0$\;
  $num \leftarrow 1$\;
  $norm \leftarrow 1$\;
  \While{$norm > \delta$ and $num < N_c$}{
    Update $u_k$ and $\mathbfcal{U}^k$ from $\sigma_k$ with CEM\;
    $\eta_b \leftarrow \frac{\|\mathbfcal{U}^t-\mathbfcal{U}^k\|_{l^2}}{\|\mathbfcal{U}^t\|_{l^2}}$\;
    \If{$\eta_b < \eta^b_0$}{
      $u_{b}\leftarrow u_k(x) \text{ for } x\in\Omega \text{ and } dist(x,\partial\Omega) = d$\;
      $\sigma \leftarrow\sigma_k$\;
      break\;
    }
    Update $\sigma_k$ with LM-SCEM\;
    $norm \leftarrow \|\tau_k\|_{L^2}$\;
    $num \leftarrow num+1$\;
  }
  Reconstruct $\mathcal{E}^\delta$ from $u_{b}$ in domain $\Omega'=\{x\ |\ x\in\Omega,\ dist(x, \partial\Omega)>d\}$\;
  $\sigma_k\leftarrow\sigma$\;
  $num \leftarrow 1$\;
  $norm \leftarrow 1$\;
  \While{$norm > \delta$ and $num < N_d$}{
    Update $u_k$ from $\sigma_k$ with DCM\;
    Update $\sigma_k$ with LM-DCM\;
    $norm \leftarrow \|\tau_k\|_{L^2}$\;
    $num \leftarrow num+1$\;
  }
  
  \caption{The LM-SCEM-DCM algorithm for reconstructing the conductivity map of a domain $\Omega$ from a single measurement of power density.}
  \label{algo_lm_cem_dcm}
\end{algorithm}

\section{Numerical investigation}\label{sec:3}
\subsection{Phantom preparation and numerical setup}
The variational forms of the linear systems defined by equations \eqref{eq:final_lma_final} in Section \ref{sec:numeri_imple} can easily be obtained through integration by parts. They are solved with a mixed finite element method \cite{Boffi2013} which is implemented  using FEniCS \cite{AlnaesBlechta2015a}. To illustrate the stability and accuracy of the presented approaches we the focus on numerical examples in a  2-dimensional (2D) problem. Two phantoms are considered.

\begin{figure}[ht]
  \centering
  \begin{subfigure}[b]{0.45\linewidth}
    \centering
    \includegraphics[height=0.5\linewidth]{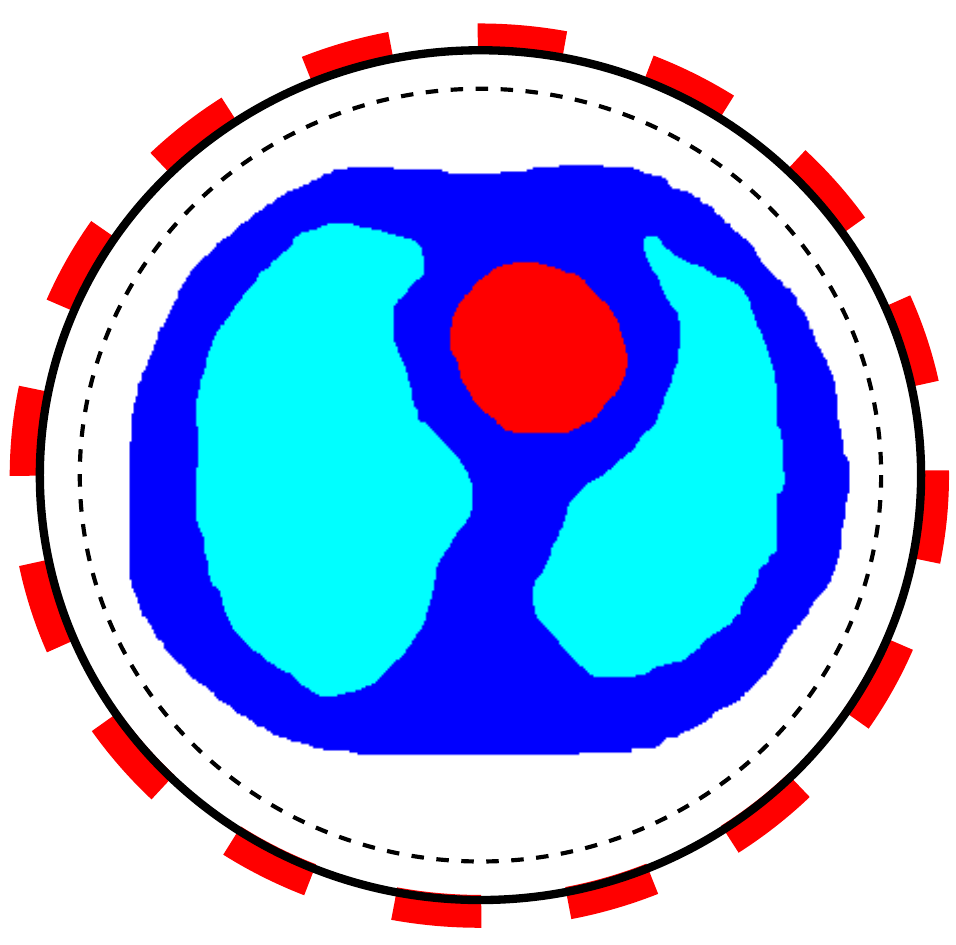}
    \caption{Heart-lung model}
    \label{2d_heart_lung}
  \end{subfigure}
  \begin{subfigure}[b]{0.45\linewidth}
    \centering
    \includegraphics[height=0.5\linewidth]{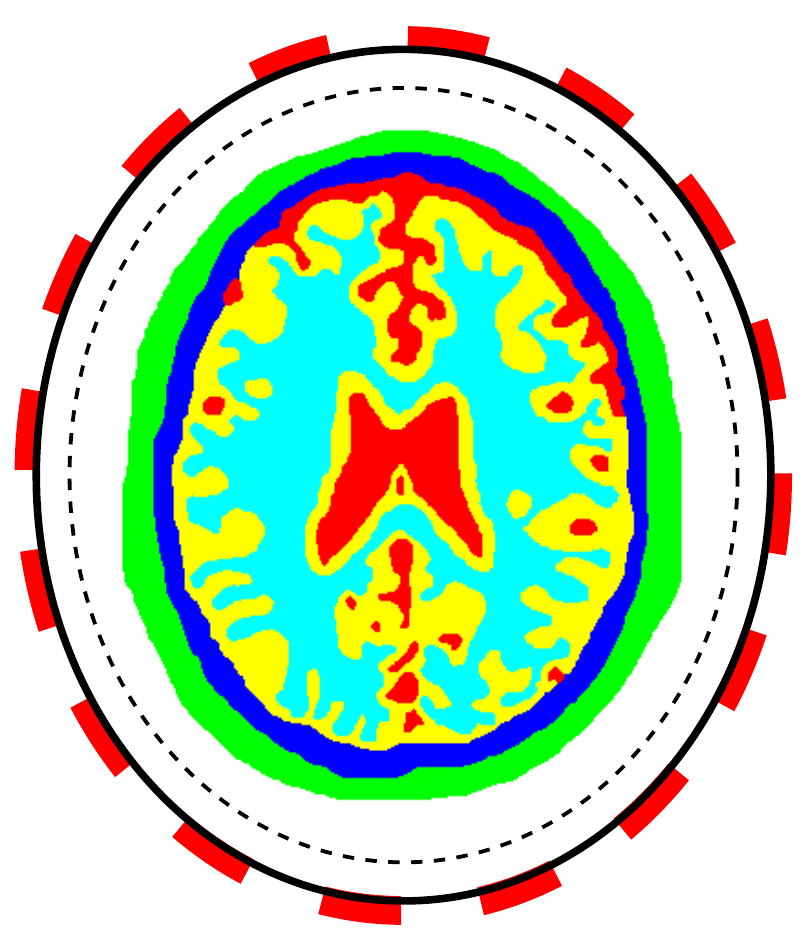}
    \caption{Human brain model}
    \label{2d_brain}
  \end{subfigure}
  \caption{The 2D (a) heart-lung model and (b) human brain model embedded in a background material with electrodes (red squares) attached to the boundary (the solid black line). Different regions are marked with different colors. In (a), there are heart (red, $\sigma=\SI{0.7}{S/m}$), lung (cyan, $\sigma=\SI{0.26}{S/m}$), soft-tissues (blue, $\sigma=\SI{0.33}{S/m}$), and background material (white, $\sigma=\SI{0.22}{S/m}$). In (b), there are scalp (green, $\sigma=\SI{0.5232}{S/m}$), skull (blue, $\sigma=\SI{0.2923}{S/m}$), cerebro-spinal fluid (red, $\sigma=\SI{2.1143}{S/m}$), gray matter (yellow, $\sigma=\SI{0.5595}{S/m}$), white matter (cyan, $\sigma=\SI{0.3240}{S/m}$) and the background mateiral (white, $\sigma=\SI{0.4}{S/m}$). All electrodes are uniformly distributed with the same corresponding central section angle. It is assumed that the electrical conductivity $\sigma$ in the region close to boundary (between solid and dashed black lines) is known.}
  \label{phantoms}
\end{figure}

\begin{figure}[ht]
  \centering
  \begin{subfigure}[b]{0.45\linewidth}
    \centering
    \includegraphics[height=0.5\linewidth]{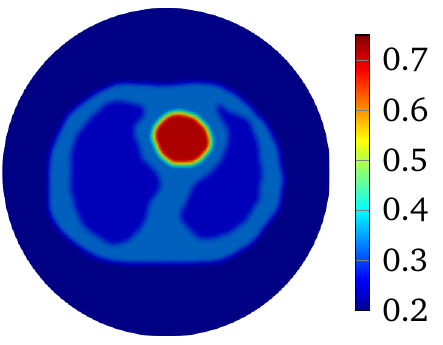}
    \caption{Heart-lung model}
    \label{2d_heartlung_sigma}
  \end{subfigure}
  \begin{subfigure}[b]{0.45\linewidth}
    \centering
    \includegraphics[height=0.5\linewidth]{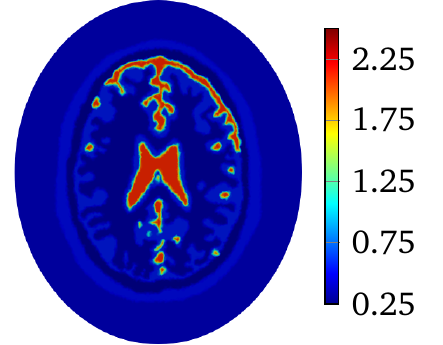}
    \caption{Human-brain model}
    \label{2d_brain_sigma}
  \end{subfigure}
  \caption{The true distribution of $\sigma$ (a) for the 2D heart-lung model shown in Figure \ref{2d_heart_lung} and (b) for the 2D human brain model shown in Figure \ref{2d_brain}}
  \label{true_sigma}
\end{figure}

The first example is a heart-lung model \cite{Zlochiver2003}, see Figure \ref{2d_heart_lung}. The considered three tissues are heart (red, $\sigma=\SI{0.7}{S/m}$), lung (cyan, $\sigma=\SI{0.26}{S/m}$), and soft-tissues (blue, $\sigma=\SI{0.33}{S/m}$). The model is placed into a circular region with a background material (white, $\sigma=\SI{0.22}{S/m}$) and a radius $r=\SI{25}{cm}$. 
The second example is the human brain model  shown in Figure \ref{2d_brain}. The considered tissues in this model include scalp (green, $\sigma=\SI{0.5232}{S/m}$), skull (blue, $\sigma=\SI{0.2923}{S/m}$), cerebro-spinal fluid (red, $\sigma=\SI{2.1143}{S/m}$), gray matter (yellow, $\sigma=\SI{0.5595}{S/m}$) and white matter (cyan, $\sigma=\SI{0.3240}{S/m}$). Refer to \cite{Andreuccetti1997} for conductivities of different tissues. The shape of this model is close to an ellipse whose semi-major and semi-minor axes are \SI{6}{cm} and \SI{7}{cm}. The model is placed in an ellipse region with a background material (white, $\sigma=\SI{0.4}{S/m}$). The semi-minor and semi-major axes of the region are \SI{8}{cm} and \SI{9}{cm}, respectively. The conductivity maps for the phantoms are piece-wise constant functions which will be mollified with
\begin{equation}
  \label{mollifier}
  \eta_\epsilon(x,y) = C\exp\left(\frac{\epsilon^2}{ (x^2 + y^2) - \epsilon^2}\right)
\end{equation}
to produce $\sigma\in\HO{s}$ with $s=2$ for 2D problem. The constant $C>0$ is selected so that $\int_{\R^2}\eta_\epsilon = 1$. The value of $\epsilon$ is \SI{1}{cm} and \SI{0.06}{cm} for heart-lung model and human-brain model, respectively. The true smoothed distributions of $\sigma$ for the two models are shown in Figure \ref{true_sigma}.

To work with CEM and LM-SCEM, 16 electrodes (red rectangles shown in Figure \ref{phantoms}) are uniformly attached on the boundary (solid black lines). The section occupied by each electrode has the same central angle in both the heart-lung and human-brain models. The following computations assume that the conductivity in a small region close to boundary (between solid and dashed black lines) is known. The distance between the solid and dashed lines is given by $\delta_d$. This known region helps to improve the convergence of the algorithm. Three current patterns based on Fourier basis functions are used in the computations, which are $I^{(n)}_l=\cos(n \theta_l)$ for $n=1,2,3$, and $\theta_l=2\pi l/L$. The regularization parameter $\alpha_k$ is chosen to decrease exponentially, and $\alpha_k={\alpha_0}/{a^k}$ with $a > 1$. In what follows, a relatively large value is given to $\alpha_0$, and a value close to 1 is given to $a$ for a slow decreasing of $\alpha_k$ to ensure the convergence of the iterations.

\subsection{Performance of LM-SCEM}
Numerical experiments on heart-lung model are carried out here to investigate the performance of LM-SCEM. A mesh of the circular domain with 77101 triangles is used for the reconstruction. The power density for each current pattern is simulated with CEM, and Gaussian white noise is added to avoid an inverse crime. Here, the noise level is measure with signal to noise ratio (SNR) $\text{SNR}=20\log_{10}\frac{\|\mathcal{E}_p\|_{L^2}}{\|\mathcal{N}\|_{L^2}}$, where $\mathcal{N}$ is a Gaussian white noise distribution. The Levenberg-Marquardt iteration is stopped when $\|\tau\|_{L^2} < \num{1e-4}$ or total number of iterations greater than 15. These values were chosen to balance the quality in the reconstructions versus the computational speed. The reconstruction with different current patterns and different level of noise are considered to check the convergence of $\sigma_k$ and the stability of LM-SCEM. The parameters for the reconstruction are given as $\alpha_0=50$, $a=1.2$, and $\beta=\num{1.2e-3}$. These parameters are chosen to ensure that the reconstruction with current pattern $I^{(2)}$ converges. Other measurements are taken into the reconstruction without changing the parameters.

\begin{figure}[ht]
  \centering
  \begin{subfigure}[b]{0.32\linewidth}
    \centering
    \includegraphics[width=\linewidth]{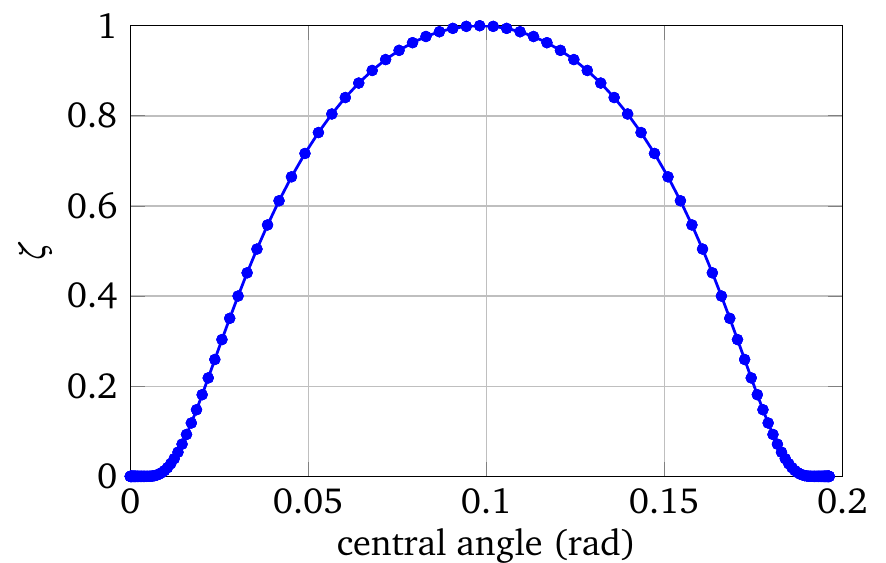}
    \caption{The distribution of $\zeta$}
    \label{zeta}
  \end{subfigure}
  \hfill
  \begin{subfigure}[b]{0.32\linewidth}
    \centering
    \includegraphics[width=\linewidth]{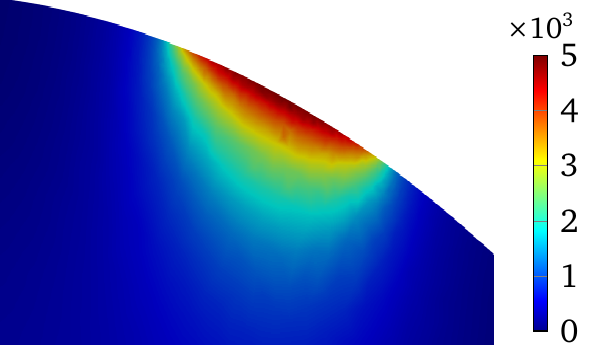}
    \caption{$\E_p(\sigma)$ from SCEM}
    \label{powerdensity_scem}
  \end{subfigure}
  \begin{subfigure}[b]{0.32\linewidth}
    \centering
    \includegraphics[width=\linewidth]{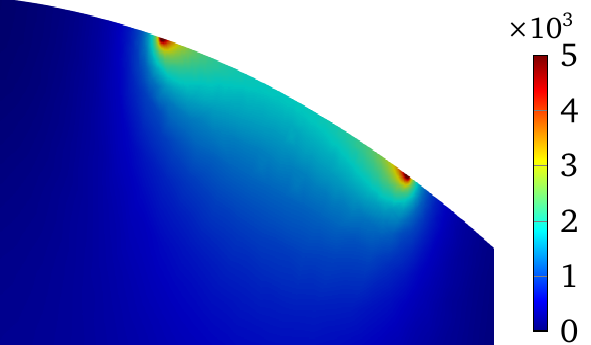}
    \caption{$\E_p(\sigma)$ from CEM}
    \label{power_density_cem}
  \end{subfigure}
  \caption{(a) The distribution of $\zeta$ on each electrode. (b) The power density $\E_p(\sigma)$ calculated from SCEM with the distribution of $\zeta$ given in (a). (c) The power density $\E_p(\sigma)$ calculated from CEM with $z_l = \frac{2}{\max\zeta}$}
  \label{CEM_Pure1}
\end{figure}

To simulate the power density $\E(\sigma)$, SCEM is used to yield better regularity of the electrical potential. The electrode conductance $\zeta_l$ of $e_l$ is chosen as
\begin{equation}
  \zeta_l(x) = \frac{1}{\epsilon^2}\exp\left(\frac{\epsilon^2}{x^2-\epsilon^2}\right),\ -\frac{l_e}{2} < x < \frac{l_e}{2},\ \epsilon=\frac{l_e}{2},\ l=1,2,\ldots,L.
\end{equation}
Here, $l_e$ is the length of the electrode, which is proportional to the central angle corresponding to the electrode. $\zeta_l(x)$ is then scaled to have the required maximum value. The distribution of $\zeta_l(x)$ for the calculation in this paper is given in Figure \ref{zeta} with a maximum value 1. The power density simulated with SCEM is given in Figure \ref{powerdensity_scem} for a region near one electrode. Here, the current pattern is $I^{(2)}$. 
A smooth distribution of $\E(\sigma)$ near the edge of the electrode is observed. For the same region, the power density simulated from CEM is also given in Figure \ref{power_density_cem} with $z_l=2.0/\max\zeta(x) = \SI{2}{S/m}$ for $l = 1,2,\ldots,L$. Two singularity points are observed at the edges of the electrode. These singularities will cause instability in the reconstruction, and therefore SCEM is used in LM-SCEM.

\begin{figure}[ht]
  \centering
  \begin{subfigure}[b]{0.45\linewidth}
    \centering
    \includegraphics[width=\linewidth]{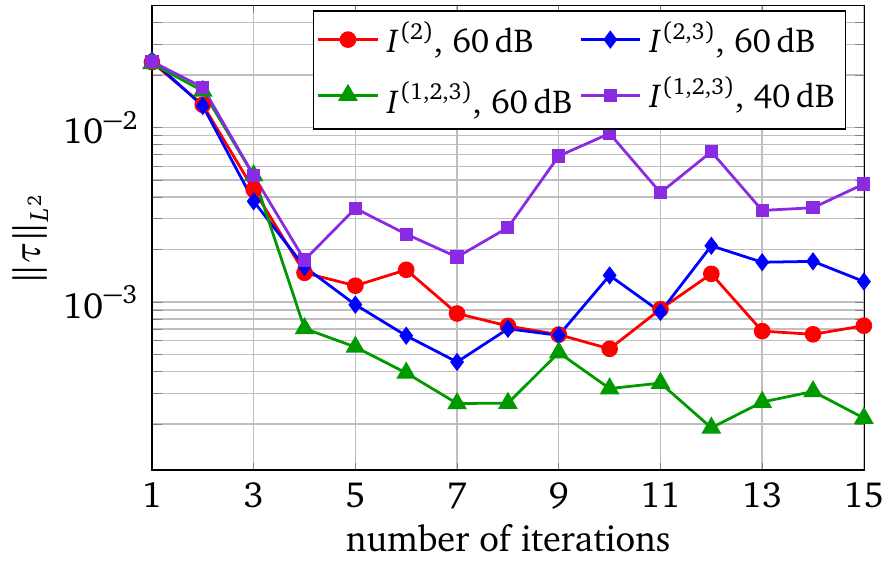}
    \caption{}
    \label{con_lm_cem}
  \end{subfigure}
  \hfill
  \begin{subfigure}[b]{0.45\linewidth}
    \centering
    \includegraphics[width=\linewidth]{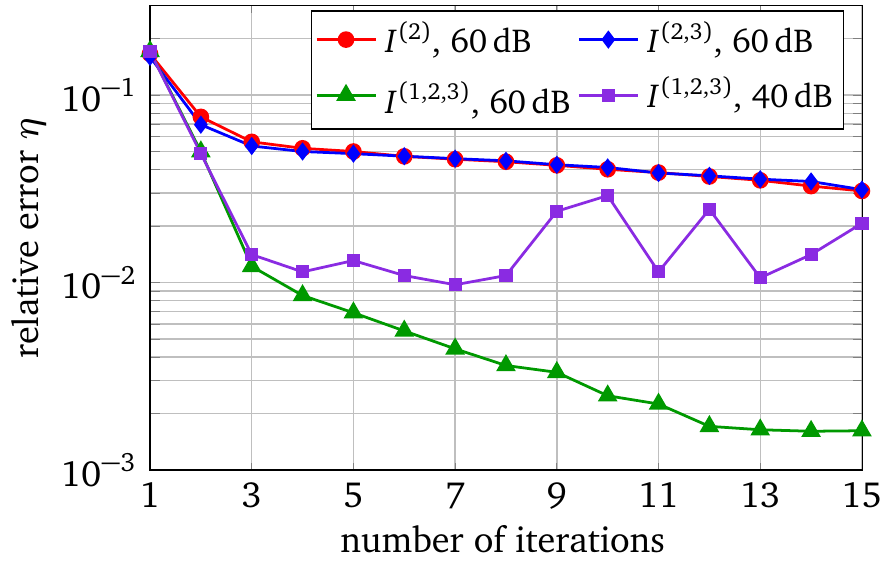}
    \caption{}
    \label{relative_error_cem}
  \end{subfigure}
  \caption{(a) The variation of $\|\tau_k\|_{L^2}$. (b) The relative error $\eta$ of the reconstruction. With \SI{60}{dB} noise, a uniform convergence is observed. With \SI{40}{dB} noise, the reconstruction almost uniformly converges to $\eta = 0.973\%$ with 7 steps.}
  \label{CEM_Pure2}
\end{figure}

\begin{figure}[ht]
  \centering
  \begin{subfigure}[b]{0.225\linewidth}
    \centering
    \includegraphics[width=\linewidth]{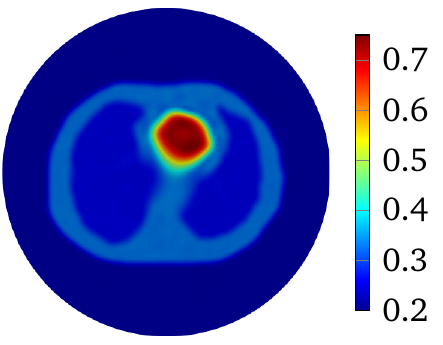}
    \caption{$I^{(2)}$, SNR = \SI{60}{dB}}
    \label{I260db}
  \end{subfigure}
  \hfill
  \begin{subfigure}[b]{0.225\linewidth}
    \centering
    \includegraphics[width=\linewidth]{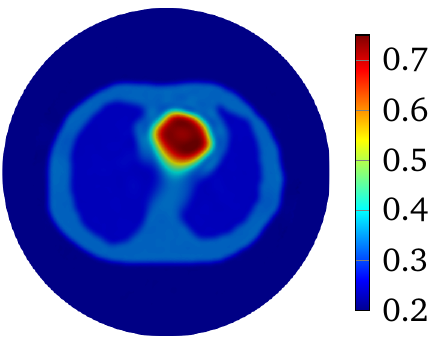}
    \caption{$I^{(2,3)}$, SNR = \SI{60}{dB}}
    \label{I2360db}
  \end{subfigure}
    \hfill
  \begin{subfigure}[b]{0.225\linewidth}
    \centering
    \includegraphics[width=\linewidth]{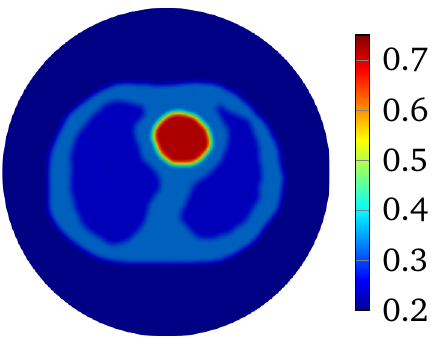}
    \caption{$I^{(1,2,3)}$, SNR = \SI{60}{dB}}
    \label{I12360db}
  \end{subfigure}
    \hfill
  \begin{subfigure}[b]{0.225\linewidth}
    \centering
    \includegraphics[width=\linewidth]{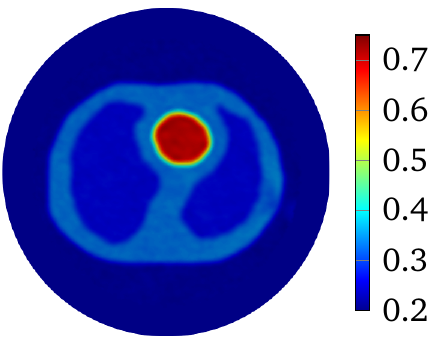}
    \caption{$I^{(1,2,3)}$, SNR = \SI{40}{dB}}
    \label{I12340db}
  \end{subfigure}
  \caption{The conductivity map reconstructed from LM-SCEM with current pattern (a) $I^{(2)}$ with \SI{60}{dB} noise, (b) $I^{(2,3)}$ with \SI{60}{dB} noise, (c) $I^{(1,2,3)}$ with \SI{60}{dB} noise, and (d) $I^{(1,2,3)}$ with \SI{40}{dB} noise.}
  \label{CEM_Pure3}
\end{figure}

The LM-SCEM algorithm is used to reconstruct the distribution of $\sigma$. The initial guess is given as $\sigma_0 = \SI{0.22}{S/m}$ which is the value of the background tissue. The conductivity in the region $x^2+y^2>(r-\delta_d)^2$ for $(x,y)\in\Omega$ is supposed to be known. The values of $\tau_k$ are truncated to only update $x^2+y^2\le (r-\delta_d)^2$ for $(x,y)\in\Omega$ with $\delta_d=\SI{4.5}{cm}$. The discontinuity caused by this truncation can be removed by applying \eqref{mollifier} properly (either by mollification or simply by replacing the discontinuous values), but it does not cause any numerical problems since $\|\tau\|_{k}$ is small, so no special treatment was done in the following computation. The values of $\|\tau\|_{L^2}$ for the first 15 iterations are shown in Figure \ref{con_lm_cem}. The relative error $\eta$ is also given in Figure \ref{relative_error_cem}. With \SI{60}{dB} noise, the reconstruction with $I^{(2)}$ uniformly converges to $\eta = 3.08\%$ with 15 iterations. The reconstructed $\sigma$ is shown in Figure \ref{I260db}. To achieve a level of $\eta = 0.1\%$, it takes more than 40 iterations. A reconstruction with $I^{(2)}$ and $I^{(3)}$ is then carried out, but a similar speed of convergence and relative error is observed, result is in Figure \ref{I2360db}. When the current pattern $I^{(1)}$ is further considered into the reconstruction, an obvious improvement of convergence is seen, and a relative error level $\eta = 0.162\%$ is achieved with 14 iterations, the conductivity map is shown in Figure \ref{I12360db}. Therefore, the convergence of LM-SCEM depends not only on the regularization parameter $\alpha_k$ and the scaling parameter $\beta$, but also on the current patterns for the measurements. Since $\text{SNR}=\SI{60}{dB}$ indicates 0.1\% of Gaussian noise in the simulated $\E_p$, the reconstructed result with $\eta = 0.162\%$ is already good, and more iterations will not improve the result. To verify, a reconstruction with $I^{(1)}$, $I^{(2)}$ and $I^{(3)}$ is carried out with $\text{SNR}=\SI{40}{dB}$ (1\% noise). Relative errors are shown in Figure \ref{relative_error_cem}, the reconstruction converges to $\eta = 0.973\%$ within 7 iterations, and more iterations did not bring any improvements. The reconstructed result is given in Figure \ref{I12340db}.

\subsection{Performance of the mixed reconstruction approach}
\begin{figure}[ht]
  \centering
  \begin{subfigure}[b]{0.45\linewidth}
    \centering
    \includegraphics[width=\linewidth]{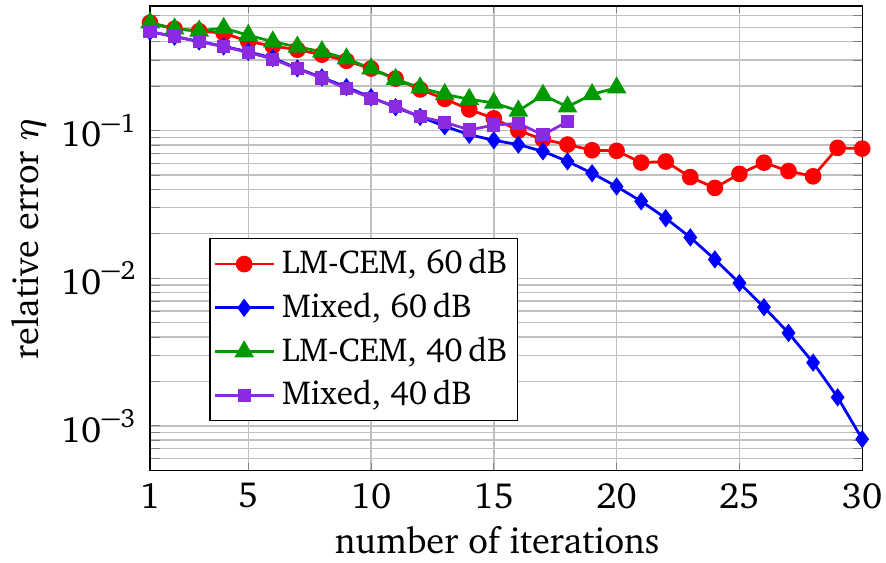}
    \caption{}
    \label{mixed_relative_error}
  \end{subfigure}
  \hfill
  \begin{subfigure}[b]{0.45\linewidth}
    \centering
    \includegraphics[width=\linewidth]{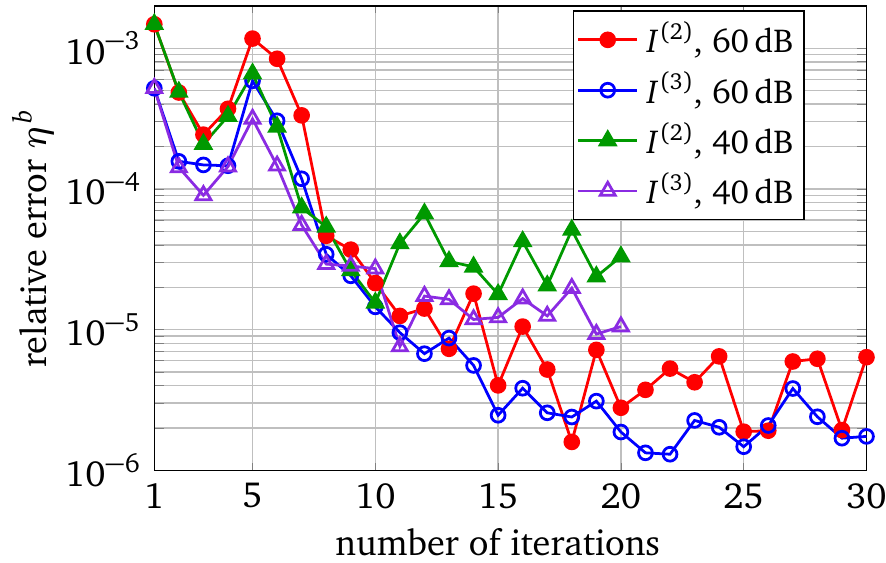}
    \caption{}
    \label{mixed_bnd_u_convergence}
  \end{subfigure}
  \caption{(a) The relative error $\eta$ of different reconstructions. (b) The convergence of the boundary potential for LM-SCEM with different level of noise.}
  \label{recon_mixed}
\end{figure}

\begin{figure}[ht]
  \centering
  \begin{subfigure}[b]{0.246\linewidth}
    \centering
    \includegraphics[width=\linewidth]{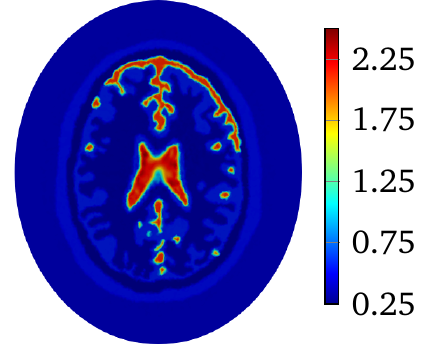}
    \caption{LM-SCEM, SNR = \SI{60}{dB}}
    \label{cem_60db_c2c3}
  \end{subfigure}
  \hfill
  \begin{subfigure}[b]{0.245\linewidth}
    \centering
    \includegraphics[width=\linewidth]{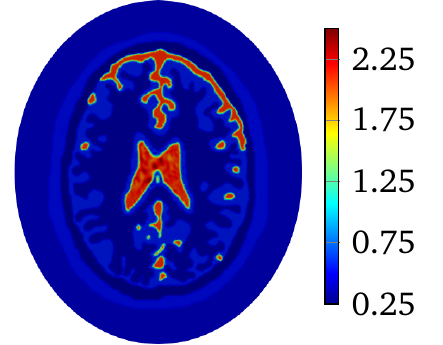}
    \caption{Mixed, SNR = \SI{60}{dB}}
    \label{mixed_60db_c2c3}
  \end{subfigure}
    \hfill
  \begin{subfigure}[b]{0.245\linewidth}
    \centering
    \includegraphics[width=\linewidth]{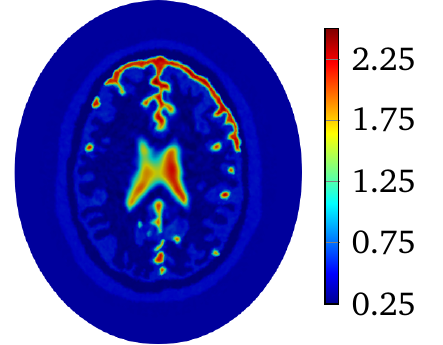}
    \caption{LM-SCEM, SNR = \SI{40}{dB}}
    \label{cem_40db_c2c3}
  \end{subfigure}
    \hfill
  \begin{subfigure}[b]{0.245\linewidth}
    \centering
    \includegraphics[width=\linewidth]{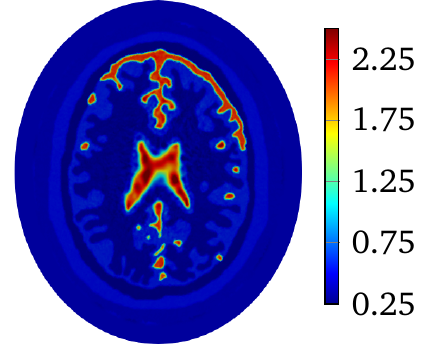}
    \caption{Mixed, SNR = \SI{40}{dB}}
    \label{mixed_40db_c2c3}
  \end{subfigure}
  \caption{The conductivity map of human-brain model reconstructed with current pattern $I^{(2)}$ and $I^{(3)}$. The reconstruction with only LM-SCEM are given in (a) and (c) for \SI{60}{dB} and \SI{40}{dB} noises. Corresponding reconstructions by mixing LM-SCEM and LM-DCM are given in (b) and (d).}
  \label{CEM_Pure4}
\end{figure}

Though LM-SCEM performs well for reconstructing the conductivity map, its efficiency of the computation decreases quickly with the increase of measurements. Since EIT is not very sensitive to the change of interior conductivity, the electrical potential should converge faster than the convergence of the conductivity. The human-brain model is  used here for numerical experiments with LM-SCEM. The elliptic domain is characterized with \SI{8}{cm} major and \SI{9}{cm} minor axes, and the domain is meshed with 36893 triangular elements. The parameters are given as $\alpha_0=150$, $a = 2.0$ and $\beta = \num{1e-3}$. The current patterns $I^{(2)}$ and $I^{(3)}$ are used in the reconstruction. With \SI{60}{dB} noise, the iteration is terminated when $\|\tau\|_{L^2}<\num{1e-4}$ or the maximum number of iterations equals 30. The reconstruction based on LM-SCEM is shown in Figure \ref{cem_60db_c2c3}. $\mathbfcal{U}^k$ of each current pattern is computed with SCEM and $\sigma_k$ for $k$-th iteration, $\eta^b$ easily follows then. The variation of the relative error $\eta$ and $\eta^b$ are shown in Figure \ref{mixed_relative_error} and Figure \ref{mixed_bnd_u_convergence}, respectively. As can be seen, $\eta$ slowly converges to $4.09\%$. This error is much larger than the input noise level, this is mainly caused by the complexity of the phantom and the high contrast of $\sigma$ among different tissues. But the potentials on the boundary for both $I^{(2)}$ and $I^{(3)}$ converge fast to a level of $\eta^b<\num{1e-4}$ in few iterations. Therefore, the boundary potential converges much faster. This property is exploited here to accelerate the computation by mixing LM-SCEM and LM-DCM, as demonstrated in Algorithm \ref{algo_lm_cem_dcm}. In this computation, the LM-SCEM is stoped when $\eta^b $ for all current patterns are smaller than $\num{1e-3}$. The LM-DCM is performed in the region $\Omega'=\{x\ |\ x\in\Omega,\ dist(x, \partial\Omega)>\SI{5}{mm}\}$ with $\delta_d = \SI{5}{mm}$. The potential on $\partial \Omega'$ is computed with SCEM and the reconstructed $\sigma$ from LM-SCEM. The power density in $\Omega'$ can be reconstructed with the method introduced by Ammari {\em et al} \cite{Ammari2008}. However, it requires the knowledge on the deformation caused by the ultrasonic waves, therefore, we compute it with DCM instead. Noise with SNR$=\SI{60}{dB}$ is added, and LM-DCM is used for the reconstruction. The relative error $\eta$ is given in Figure \ref{mixed_relative_error}. The conductivity map is reconstructed with $\eta=\num{8.13e-4}$ in 30 iterations, and the result is given in Figure \ref{mixed_60db_c2c3}. Here, the time required for 30 LM-DCM iterations is about 20 minutes which is approximately the time needed for one LM-SCEM iteration. So the reconstruction efficiency is greatly improved, and better results are obtained. A similar computation with \SI{40}{dB} noise is further considered here. As seen in Figure \ref{mixed_bnd_u_convergence}, increasing noise does not influence much the convergence of the boundary potential, therefore, this mixed approach can be a good way to remove noise from the measured power density. With 40dB noise in the reconstructed power density in $\Omega'$, the distribution of $\sigma$ obtained with LM-DCM is shown in Figure \ref{mixed_40db_c2c3}. Comparing it to the results obtained with LM-SCEM, as shown in Figure \ref{cem_40db_c2c3}, a better noise tolerance is observed in LM-DCM.

\section{Conclusion}\label{sec:4}
The work first developed a computational approach for AET by incorporating CEM into the Levenberg-Marquardt algorithm. Since the regularity of the power density obtained with traditional CEM is limited because of the Robin-type discontinues boundary conditions, a recently proposed smoothed CEM is used in this paper. Numerical investigation shows that this iterative method can stably reconstruct the conductivity map of complex phantoms, and a good accuracy can be achieved even with a certain level of noise though this also depends on the current patterns used in the measurements.

However, the reconstruction algorithm becomes inefficient quickly when the number of measurements is increased. Since EIT is not very sensitive to the changes of the internal conductivity, the boundary potential converges to its true value much faster than the convergence of the conductivity. This fact is then used to build a mixed computational approach. LM-SCEM is used to reconstruct the boundary potentials in a few iterations, and LM-DCM is applied to reconstruct the conductivity distribution based on the obtained boundary potential. It is observed in the given example that the time for one LM-SCEM iteration is enough for the whole calculation of LM-DCM, so the reconstruction efficiency is greatly improved with the mixed strategy.

Here, reconstruction with LM-DCM requires a first step to compute the power density from the obtained boundary potential from LM-SCEM. This additional step can be considered as a step to smooth the electrical potential and to remove the noise in the measured power density, in addition to the better noise tolerance of LM-DCM compared to LM-SCEM, a better reconstruction can be obtained with the mixed algorithm.

Several numerical experiments are carried out with a heart-lung phantom and a complex human brain model. The performance of the presented reconstruction approaches are well demonstrated with different number of measurements. The proposed method applies also to 3-dimensional acousto-electric tomography and  anisotropic conductivy distributions; we leave the implementation in these scenarious for future work.

\section*{Acknowledgements}
  The majority of the work was done while Changyou Li was postdoc at the  Department of Applied Mathematics and Computer Science, Technical University of Denmark in 2018. We thank  the  Danish Council for Independent Research | Natural Sciences (grant 4002-00123) for financial support.



\end{document}